\documentclass[10pt]{article}
\usepackage{amsfonts}
\usepackage{amssymb}
\newcommand{\fidemo}{\vrule height4pt width3pt depth2pt}
\newcommand{\rees}[1]{\mbox{$\mathcal{R}(#1)$}}
\newcommand{\mrees}[2]{\mbox{$\mathcal{R}(#1,#2)$}}
\newcommand{\reesw}[2]{\mbox{$\mathcal{R}(#1;#2)$}}
\newcommand{\mreesw}[3]{\mbox{$\mathcal{R}(#1,#2;#3)$}}
\newcommand{\agr}[1]{\mbox{$\mathcal{G}(#1)$}}
\newcommand{\magr}[2]{\mbox{$\mathcal{G}(#1,#2)$}}
\newcommand{\agrw}[2]{\mbox{$\mathcal{G}(#1;#2)$}}
\newcommand{\magrw}[3]{\mbox{$\mathcal{G}(#1,#2;#3)$}}
\newcommand{\tor}[4]{\mbox{${\rm Tor}_{#1}^{#2}(#3,#4)$}}

\newtheorem{rteorema}{Theorem}
\newtheorem{rcorollari}[rteorema]{Corollary}

\newtheorem{teorema}{Theorem}[section]

\newtheorem{proposicio}[teorema]{Proposition}
\newtheorem{lema}[teorema]{Lemma}
\newtheorem{corollari}[teorema]{Corollary}
\newtheorem{observacio}[teorema]{Remark}
\newtheorem{exemple}[teorema]{Example}

\begin{document}

\title{Normal transversality and uniform bounds}

\author{{\sc Francesc Planas-Vilanova} \\ \\ {\footnotesize
Dept. Matem\`atica Aplicada 1. ETSEIB-UPC.  Diagonal 647, 08028
Barcelona. E-mail:planas@ma1.upc.es}  
}

\date{}

\maketitle

\section{Introduction}

Let $A$ be a commutative ring. A graded $A$-algebra $U=\oplus _{n\geq
0}U_{n}$ is a {\em standard} $A$-algebra if $U_{0}=A$ and $U=A[U_{1}]$
is generated as an $A$-algebra by the elements of $U_{1}$. A graded
$U$-module $F=\oplus _{n\geq 0}F_{n}$ is a {\em standard} $U$-module
if $F$ is generated as an $U$-module by the elements of $F_{0}$, that
is, $F_{n}=U_{n}F_{0}$ for all $n\geq 0$.  In particular,
$F_{n}=U_{1}F_{n-1}$ for all $n\geq 1$.  Given $I$, $J$, two ideals of
$A$, we consider the following standard algebras: the {\em Rees
algebra of} $I$, $\rees{I}=\oplus _{n\geq 0}I^{n}t^{n}=A[It]\subset
A[t]$, and the {\em multi-Rees algebra of} $I$ and $J$,
$\mrees{I}{J}=\oplus _{n\geq 0}(\oplus _{p+q=n}I^{p}J^{q}u^
{p}v^{q})=A[Iu,Jv]\subset A[u,v]$. Consider the {\em associated graded
ring of} $I$, $\agr{I}=\rees{I}\otimes A/I=\oplus _{n\geq
0}I^{n}/I^{n+1}$, and the {\em multi-associated graded ring of} $I$
and $J$, $\magr{I}{J}=\mrees{I}{J}\otimes A/(I+J)=\oplus _{n\geq
0}(\oplus _{p+q=n}I^{p}J^{q}/(I+J)I^{p}J^{q})$. We can always consider
the tensor product of two standard $A$-algebras $U=\oplus _{p\geq
0}U_{p}$ and $V=\oplus _{q\geq 0}V_{q}$ as an standard $A$-algebra
with the natural grading $U\otimes V=\oplus _{n\geq 0}(\oplus
_{p+q=n}U_{p}\otimes V_{q})$. If $M$ is an $A$-module, we have the
standard modules: the {\em Rees module of} $I$ {\em with respect to}
$M$, $\reesw{I}{M}=\oplus _{n\geq 0}I^{n}Mt^{n}=M[It]\subset M[t]$ (a
standard $\rees{I}$-module), and the {\em multi-Rees module of} $I$
and $J$ {\em with respect to} $M$, $\mreesw{I}{J}{M}=\oplus _{n\geq
0}(\oplus _{p+q=n}I^{p}J^ {q}Mu^{p}v^{q})=M[Iu,Jv]\subset M[u,v]$ (a
standard $\mrees{I}{J}$-module). Consider the {\em associated graded
module of} $I$ {\em with respect to} $M$,
$\agrw{I}{M}=\reesw{I}{M}\otimes A/I=\oplus _{n\geq 0}I^{n}M/I^{n+1}M$
(a standard $\agr{I}$-module), and the {\em multi-associated graded
module of} $I$ and $J$ {\em with respect to} $M$,
$\magrw{I}{J}{M}=\mreesw{I}{J}{M}\otimes A/(I+J)=\oplus _{n\geq
0}(\oplus _{p+q=n}I^{p}J^{q}M/(I+J)I^{p}J^{q}M)$ (a standard
$\mrees{I}{J}$-module). If $U$, $V$ are two standard $A$-algebras and
$F$ is a standard $U$-module and $G$ is a standard $V$-module, then
$F\otimes G=\oplus _{n\geq 0}(\oplus _{p+q=n}F_{p}\otimes G_{q})$ is a
standard $U\otimes V$-module.

Denote by $\pi :\rees{I}\otimes \reesw{J}{M}\rightarrow
\mreesw{I}{J}{M}$ and $\sigma :\mreesw{I}{J}{M}\rightarrow
\reesw{I+J}{M}$ the natural surjective graded morphisms of standard
$\rees{I}\otimes \rees{J}$-modules. Let $\varphi :\rees{I}\otimes
\reesw{J}{M}\rightarrow \reesw{I+J}{M}$ be $\sigma \circ \pi$.  Denote
by $\overline{\pi} :\agr{I}\otimes \agrw{J}{M}\rightarrow
\magrw{I}{J}{M}$ and $\overline{\sigma} :\magrw{I}{J}{M}\rightarrow
\agrw{I+J}{M}$ the tensor product of $\pi$ and $\sigma$ by $A/(I+J)$;
these are two natural surjective graded morphisms of standard
$\agr{I}\otimes \agr{J}$-modules. Let $\overline{\varphi}
:\agr{I}\otimes \agrw{J}{M}\rightarrow \agrw{I+J}{M}$ be
$\overline{\sigma}\circ \overline{\pi}$.  The first purpose of this
note is to prove the following theorem:

\begin{rteorema}\label{rtrans}
Let $A$ be a noetherian ring, $I$, $J$ two ideals of $A$ and $M$ a
finitely generated $A$-module. The following two conditions are
equivalent:
\begin{itemize}
\item[$(i)$] $\overline{\varphi}:\agr{I}\otimes \agrw{J}{M}\rightarrow
\agrw{I+J}{M}$ is an isomorphism.
\item[$(ii)$] $\tor{1}{}{A/I^{p}}{\reesw{J}{M}}=0$ and
$\tor{1}{}{A/I^{p}}{\agrw{J}{M}}=0$ for all integers $p\geq 1$.
\end{itemize}
In particular, $\agr{I}\otimes \agr{J}\simeq \agr{I+J}$ if and only if
$\tor{1}{}{A/I^{p}}{A/J^{q}}=0$ and $\tor{2}{}{A/I^{p}}{A/J^{q}}=0$
for all integers $p,q\geq 1$.
\end{rteorema}

The morphism $\overline{\varphi}$ has been studied by Hironaka
\cite{hironaka}, Grothendieck \cite{ega} and Hermann,
Ikeda and Orbanz \cite{hio}, among others, but assuming always $A$ is
normally flat along $I$ (see 21.11 in \cite{hio}). We will see how
Theorem \ref{rtrans} generalizes all this former work.

Let us now recall some definitions in order to state the second
purpose of this note. If $U$ is a standard $A$-algebra and $F$ is a
graded $U$-module, put $s(F)={\rm min}\{ r\geq 1\mid F_{n}=0 \mbox{
for all } n\geq r+1\} $, where $s(F)$ may possibly be infinite. If
$U_{+}=\oplus _{n>0}U_{n}$ and $r\geq 1$, the following three
conditions are equivalent: $F$ can be generated by elements of degree
at most $r$; $s(F/U_{+}F)\leq r$; and $F_{n}=U_{1}F_{n-1}$ for all
$n\geq r+1$. If $\varphi :G\rightarrow F$ is a surjective graded
morphism of graded $U$-modules, we denote by $E(\varphi )$ the graded
$A$-module $E(\varphi )={\rm ker}\varphi /U_{+}{\rm ker}\varphi ={\rm
ker}\varphi _{0}\oplus (\oplus _{n\geq 1}{\rm ker}\varphi
_{n}/U_{1}{\rm ker}\varphi _{n-1})=\oplus _{n\geq 0}E(\varphi )_{n}$.
If $F$ is a standard $U$-module, take ${\bf S}(U_{1})$ the symmetric
algebra of $U_{1}$, $\alpha :{\bf S}(U_{1})\rightarrow U$ the
surjective graded morphism of standard $A$-algebras induced by the
identity on $U_{1}$ and $\gamma :{\bf S}(U_{1})\otimes F_{0}\buildrel
\alpha \otimes 1\over \rightarrow U\otimes F_{0}\rightarrow F$ the
composition of $\alpha \otimes 1$ with the structural morphism. Since
$F$ is a standard $U$-module, $\gamma $ is a surjective graded
morphism of graded ${\bf S}(U_{1})$-modules. The {\em module of
effective} $n$-{\em relations} of $F$ is defined to be
$E(F)_{n}=E(\gamma )_{n}= {\rm ker}\gamma _{n}/U_{1}{\rm ker}\gamma
_{n-1}$ (for $n=0$, $E(F)_{n}=0$). Put $E(F)=\oplus _{n\geq
1}E(F)_{n}=\oplus _{n\geq 1}E(\gamma )_{n}=E(\gamma )= {\rm ker}\gamma
/{\bf S}_{+}(U_{1}){\rm ker}\gamma$. The {\em relation type} of $F$ is
defined to be ${\rm rt}(F)=s(E(F))$, that is, ${\rm rt}(F)$ is the
minimum positive integer $r\geq 1$ such that the effective
$n$-relations are zero for all $n\geq r+1$. A {\em symmetric
presentation} of a standard $U$-module $F$ is a surjective graded
morphism of standard $V$-modules $\varphi :G\rightarrow F$, with
$\varphi :G=V\otimes M\buildrel f\otimes h\over \rightarrow U\otimes
F_{0}\rightarrow F$, where $V$ is a symmetric $A$-algebra,
$f:V\rightarrow U$ is a surjective graded morphism of standard
$A$-algebras, $h:M\rightarrow F_{0}$ is an epimorphism of $A$-modules
and $U\otimes F_{0}\rightarrow F$ is the structural morphism. One can
show (see \cite{planas2}) that $E(F)_{n}=E(\varphi )_{n}$ for all
$n\geq 2$ and $s(E(F))=s(E(\varphi ))$. Thus the module of effective
$n$-relations and the relation type of a standard $U$-module are
independent of the chosen symmetric presentation. Roughly speaking,
the relation type of $F$ is the largest degree of any minimal
homogeneous system of generators of the submodule defining $F$ as a
quotient of a polynomial ring with coefficients in $F_{0}$. For an
ideal $I$ of $A$ and an $A$-module $M$, the module of effective
$n$-relations and the relation type of $I$ with repect to $M$ are
defined to be $E(I;M)_{n}=E(\reesw{I}{M})_{n}$ and ${\rm rt}(I;M)={\rm
rt}(\reesw{I}{M})$, respectively. Then:

\begin{rteorema}\label{rrt}
Let $A$ be a commutative ring, $U$ and $V$ two standard $A$-algebras,
$F$ a standard $U$-module and $G$ a standard $V$-module. Then
$U\otimes V$ is a standard $A$-algebra, $F\otimes G$ is a standard
$U\otimes V$-module and ${\rm rt}(F\otimes G)\leq {\rm
max}({\rm rt}(F),{\rm rt}(G))$.
\end{rteorema}

As a consequence of Theorems \ref{rtrans} and \ref{rrt}, one deduces
the existence of an uniform bound for the relation type of all maximal
ideals of an excellent ring.

\begin{rteorema}\label{rboundmax}
Let $A$ be an excellent (or $J-2$) ring and let $M$ be a finitely
generated $A$-module.  Then there exists an integer $s\geq 1$ such
that, for all maximal ideals $\mathfrak{m}$ of $A$, the relation type
of $\mathfrak{m}$ with respect to $M$ satisfies ${\rm
rt}(\mathfrak{m};M)\leq s$.
\end{rteorema}

In fact, Theorem \ref{rboundmax} could also been deduced from the
proof of Theorem 4 of Trivedi in \cite{trivedi}. Finally, and using
Theorem 2 of \cite{planas2}, one can recover the following result of
Duncan and O'Carroll.

\begin{rcorollari}{\rm \cite{do}}
Let $A$ be an excellent (or $J-2$) ring and let $N\subseteq M$ be
two finitely generated $A$-modules.
Then there exists an integer $s\geq 1$ such that, for all integers
$n\geq s$ and for all maximal ideals $\mathfrak{m}$ of $A$,
$\mathfrak{m}^{n}M\cap N=\mathfrak{m}^{n-s}(\mathfrak{m}^{s}M\cap N)$.
\end{rcorollari}

\section{Normal transversality}

\begin{lema}\label{tecnic}
Let $A$ be a commutative ring, $I$ an ideal of $A$, $U$ a standard
$A$-algebra, $F$ and $G$ two standard $U$-modules and $\varphi
:G\rightarrow F$ a surjective graded morphism of standard
$A$-alegbras. If $\overline{A}=A/I$, then $\overline{U}=U\otimes
\overline{A}$ is a standard $\overline{A}$-algebra,
$\overline{F}=F\otimes \overline{A}$ and $\overline{G}=G\otimes
\overline{A}$ are two standard $\overline{U}$-modules and
$\overline{\varphi}=\varphi \otimes 1_{\overline{A}}:
\overline{G}\rightarrow \overline{F}$ is a surjective graded morphism
of standard $\overline{U}$-modules. Moreover,
$s(E(\overline{\varphi}))\leq s(E(\varphi ))$.
\end{lema}

\noindent {\em Proof}. Consider the following commutative diagram of
exact rows:

\begin{picture}(330,85)(-5,0)

\put(100,60){\makebox(0,0){\mbox{\footnotesize
$U_{1}\otimes {\rm ker}\varphi _{n-1}$}}}
\put(200,60){\makebox(0,0){\mbox{\footnotesize
$U_{1}\otimes G_{n-1}$}}}
\put(300,60){\makebox(0,0){\mbox{\footnotesize
$U_{1}\otimes F_{n-1}$}}}
\put(360,60){\makebox(0,0){$0$}}

\put(135,60){\vector(1,0){22}}
\put(240,60){\vector(1,0){25}}
\put(330,60){\vector(1,0){20}}

\put(40,20){\vector(1,0){30}}
\put(135,20){\vector(1,0){30}}
\put(235,20){\vector(1,0){30}}
\put(325,20){\vector(1,0){25}}

\put(30,20){\makebox(0,0){$0$}}
\put(100,20){\makebox(0,0){{\footnotesize ${\rm ker}\varphi _{n}$}}}
\put(200,20){\makebox(0,0){{\footnotesize $G_{n}$}}}
\put(300,20){\makebox(0,0){{\footnotesize $F_{n}$}}}
\put(360,20){\makebox(0,0){$0$}}
\put(365,17){\makebox(0,0){.}}

\put(100,50){\vector(0,-1){20}}
\put(200,50){\vector(0,-1){20}}
\put(200,50){\vector(0,-1){16}}
\put(300,50){\vector(0,-1){20}}
\put(300,50){\vector(0,-1){16}}
\put(210,40){\makebox(0,0){\mbox{\footnotesize $\partial _{n}^{G}$}}}
\put(310,40){\makebox(0,0){\mbox{\footnotesize $\partial _{n}^{F}$}}}
\put(250,70){\makebox(0,0){\mbox{\footnotesize 
$1\otimes \varphi _{n-1}$}}}
\put(250,10){\makebox(0,0){\mbox{\footnotesize $\varphi _{n}$}}}
\end{picture}

\noindent By the snake lemma, ${\rm ker}\partial _{n}^{G}\rightarrow
{\rm ker}\partial _{n}^{F}\rightarrow E(\varphi )_{n}\rightarrow 0$ is
an exact sequence of $A$-modules. If we tensor this sequence by
$\overline{A}$, then $({\rm ker}\partial _{n}^{G})\otimes
\overline{A}\rightarrow ({\rm ker}\partial _{n}^{F})\otimes
\overline{A}\rightarrow E(\varphi )_{n}\otimes \overline{A}\rightarrow
0$ is an exact sequence of $\overline{A}$-modules. On the other hand,
we have the following commutative diagram of exact rows:

\begin{picture}(330,85)(-5,0)

\put(100,60){\makebox(0,0){\mbox{\footnotesize
$\overline{U}_{1}\otimes {\rm ker}\overline{\varphi}_{n-1}$}}}
\put(200,60){\makebox(0,0){\mbox{\footnotesize
$\overline{U}_{1}\otimes \overline{G}_{n-1}$}}}
\put(300,60){\makebox(0,0){\mbox{\footnotesize
$\overline{U}_{1}\otimes \overline{F}_{n-1}$}}}
\put(360,60){\makebox(0,0){$0$}}

\put(135,60){\vector(1,0){22}}
\put(240,60){\vector(1,0){25}}
\put(330,60){\vector(1,0){20}}

\put(40,20){\vector(1,0){30}}
\put(135,20){\vector(1,0){30}}
\put(235,20){\vector(1,0){30}}
\put(325,20){\vector(1,0){25}}

\put(30,20){\makebox(0,0){$0$}}
\put(100,20){\makebox(0,0){{\footnotesize 
${\rm ker}\overline{\varphi}_{n}$}}}
\put(200,20){\makebox(0,0){{\footnotesize $\overline{G}_{n}$}}}
\put(300,20){\makebox(0,0){{\footnotesize $\overline{F}_{n}$}}}
\put(360,20){\makebox(0,0){$0$}}
\put(365,17){\makebox(0,0){.}}

\put(100,50){\vector(0,-1){20}}
\put(200,50){\vector(0,-1){20}}
\put(200,50){\vector(0,-1){16}}
\put(300,50){\vector(0,-1){20}}
\put(300,50){\vector(0,-1){16}}
\put(210,40){\makebox(0,0){\mbox{\footnotesize 
$\partial _{n}^{\overline{G}}$}}}
\put(310,40){\makebox(0,0){\mbox{\footnotesize 
$\partial _{n}^{\overline{F}}$}}}
\put(250,70){\makebox(0,0){\mbox{\footnotesize 
$1\otimes \overline{\varphi}_{n-1}$}}}
\put(250,10){\makebox(0,0){\mbox{\footnotesize $\overline{\varphi}_{n}$}}}
\end{picture}

\noindent By the snake lemma, ${\rm ker}\partial
_{n}^{\overline{G}}\rightarrow {\rm ker}\partial
_{n}^{\overline{F}}\rightarrow E(\overline{\varphi})_{n}\rightarrow 0$
is an exact sequence of $\overline{A}$-modules. In order to see the
relationship between ${\rm ker}\partial _{n}^{F}$ and ${\rm
ker}\partial _{n}^{\overline{F}}$, tensor by $\overline{A}$ the exact
sequence of $A$-modules $0\rightarrow {\rm ker}\partial
_{n}^{F}\rightarrow U_{1}\otimes F_{n-1}\buildrel \partial
_{n}^{F}\over \rightarrow F_{n}\rightarrow 0$ and consider the
commutative diagram of exact rows:

\begin{picture}(330,85)(-5,0)

\put(100,60){\makebox(0,0){\mbox{\footnotesize
$({\rm ker}\partial _{n}^{F})\otimes \overline{A}$}}}
\put(200,60){\makebox(0,0){\mbox{\footnotesize
$(U_{1}\otimes F_{n-1})\otimes \overline{A}$}}}
\put(300,60){\makebox(0,0){\mbox{\footnotesize
$F_{n}\otimes \overline{A}$}}}
\put(360,60){\makebox(0,0){$0$}}

\put(135,60){\vector(1,0){22}}
\put(240,60){\vector(1,0){25}}
\put(330,60){\vector(1,0){20}}

\put(40,20){\vector(1,0){30}}
\put(135,20){\vector(1,0){30}}
\put(235,20){\vector(1,0){30}}
\put(325,20){\vector(1,0){25}}

\put(30,20){\makebox(0,0){$0$}}
\put(100,20){\makebox(0,0){{\footnotesize 
${\rm ker}\partial _{n}^{\overline{F}}$}}}
\put(200,20){\makebox(0,0){{\footnotesize 
$\overline{U}_{1}\otimes \overline{F}_{n-1}$}}}
\put(300,20){\makebox(0,0){{\footnotesize $\overline{F}_{n}$}}}
\put(360,20){\makebox(0,0){$0$}}
\put(365,17){\makebox(0,0){.}}

\put(200,50){\vector(0,-1){20}}
\put(300,50){\vector(0,-1){20}}
\put(210,40){\makebox(0,0){\mbox{\footnotesize $\simeq $}}}
\put(310,40){\makebox(0,0){\mbox{\footnotesize $\simeq $}}}
\put(250,70){\makebox(0,0){\mbox{\footnotesize 
$\partial _{n}^{F}\otimes 1$}}}
\put(250,10){\makebox(0,0){\mbox{\footnotesize $\partial
_{n}^{\overline{F}}$}}}
\end{picture}

\noindent It induces an epimorphism of $\overline{A}$-modules $({\rm
ker}\partial _{n}^{F})\otimes \overline{A}\rightarrow {\rm
ker}\partial _{n}^{\overline{F}}$. Analogously, there exists an
epimorphism of $\overline{A}$-modules $({\rm ker}\partial
_{n}^{G})\otimes \overline{A}\rightarrow {\rm ker}\partial
_{n}^{\overline{G}}$. Both epimorphims make commutative the following
diagram of exact rows:

\begin{picture}(330,85)(-5,0)

\put(100,60){\makebox(0,0){\mbox{\footnotesize
$({\rm ker}\partial _{n}^{G})\otimes \overline{A}$}}}
\put(200,60){\makebox(0,0){\mbox{\footnotesize
$({\rm ker}\partial _{n}^{F})\otimes \overline{A}$}}}
\put(300,60){\makebox(0,0){\mbox{\footnotesize
$E(\varphi )_{n}\otimes \overline{A}$}}}
\put(360,60){\makebox(0,0){$0$}}

\put(135,60){\vector(1,0){22}}
\put(240,60){\vector(1,0){25}}
\put(330,60){\vector(1,0){20}}

\put(135,20){\vector(1,0){30}}
\put(235,20){\vector(1,0){30}}
\put(325,20){\vector(1,0){25}}

\put(100,20){\makebox(0,0){{\footnotesize 
${\rm ker}\partial _{n}^{\overline{G}}$}}}
\put(200,20){\makebox(0,0){{\footnotesize 
${\rm ker}\partial _{n}^{\overline{F}}$}}}
\put(300,20){\makebox(0,0){{\footnotesize 
$E(\overline{\varphi})_{n}$}}}
\put(360,20){\makebox(0,0){$0$}}
\put(365,17){\makebox(0,0){,}}

\put(100,50){\vector(0,-1){20}}
\put(100,50){\vector(0,-1){16}}
\put(200,50){\vector(0,-1){20}}
\put(200,50){\vector(0,-1){16}}
\end{picture}

\noindent from where we deduce an epimorphism $E(\varphi )_{n}\otimes
\overline{A}\rightarrow E(\overline{\varphi})_{n}$. In particular,
$s(E(\overline{\varphi}))\leq s(E(\varphi ))$. \fidemo

\medskip

\begin{lema}\label{necesig}
Let $A$ be a commutative ring, $I$, $J$ two ideals of $A$ and $M$ an
$A$-module. Consider $\sigma :\mreesw{I}{J}{M}\rightarrow \reesw{I+J}{M}$
and $\overline{\sigma}=\magrw{I}{J}{M}\rightarrow \agrw{I+J}{M}$. Then
\begin{itemize}
\item[$(a)$] ${\rm ker}(\sigma _{1})\simeq IM\cap JM$.
\item[$(b)$] ${\rm ker}(\overline{\sigma}_{1})=0$ if and only if
$IM\cap JM\subset I(I+J)M\cap (I+J)JM$.
\item[$(c)$] If $I^{p}M\cap J^{q}M=I^{p}J^{q}M$ for all integers
$p,q\geq 1$, then $s(E(\sigma ))=1$ and $\overline{\sigma}$ is an
isomorphism.
\end{itemize}
\end{lema}

\noindent {\em Proof}. Consider $0\rightarrow IM\cap JM\buildrel \rho
\over \rightarrow IM\oplus JM\buildrel \sigma _{1} \over \rightarrow
(I+J)M\rightarrow 0$ where $\rho (a)=(a,-a)$ and $\sigma
_{1}(a,b)=a+b$. Clearly it is an exact sequence of $A$-modules. Thus
${\rm ker}(\sigma _{1})=\rho (IM\cap JM)\simeq IM\cap JM$. If we
tensor this exact sequence by $A/(I+J)$ we get $(IM\cap JM)\otimes
A/(I+J)\buildrel \overline{\rho} \over \rightarrow (IM\oplus
JM)\otimes A/(I+J)\buildrel \overline{\sigma} _{1} \over \rightarrow
(I+J)M/(I+J)^{2}M\rightarrow 0$. Then
\begin{eqnarray*}
{\rm ker}(\overline{\sigma}_{1})={\rm im}\overline{\rho}=\{
(\overline{a},-\overline{a})\in IM/I(I+J)M\oplus JM/(I+J)JM\mid a\in
IM\cap JM\} \, .
\end{eqnarray*}
Hence ${\rm ker}(\overline{\sigma}_{1})=0$ if and only if $IM\cap
JM\subset I(I+J)M\cap (I+J)JM$. Now, let us prove $(c)$. Let $z\in
{\rm ker}\sigma _{n}\subset \mreesw{I}{J}{M}_{n}=\oplus
_{p+q=n}I^{p}J^{q}Mu^{p}v^{q}\subset M[u,v]$. Thus,
$z=a_{0}u^{n}+a_{1}u^{n-1}v+\ldots +a_{n-1}uv^{n-1}+a_{n}v^{n}$,
$a_{i}\in I^{n-i}J^{i}M$, and $0=\sigma _{n}(z)=(a_{0}+a_{1}+\ldots
+a_{n-1}+a_{n})t^{n}\in \reesw{I+J}{M}_{n}=(I+J)^{n}Mt^{n}$. So
$a_{0}+a_{1}+\ldots +a_{n-1}+a_{n}=0$. Let us denote: 

$$\left\{ \begin{array}{l}
b_{0}=a_{0}\in I^{n}M\cap JM=I^{n}JM\\
b_{1}=a_{0}+a_{1}\in I^{n-1}M\cap J^{2}M=I^{n-1}J^{2}M\mbox{ and
}a_{1}=b_{1}-b_{0}\\
b_{2}=a_{0}+a_{1}+a_{2}\in I^{n-2}M\cap J^{3}M=I^{n-2}J^{3}M\mbox{ and
}a_{2}=b_{2}-b_{1}\\ \ldots \\
b_{n-2}=a_{0}+\ldots +a_{n-2}\in I^{2}M\cap
J^{n-1}M=I^{2}J^{n-1}M\mbox{ and }a_{n-2}=b_{n-2}-b_{n-3}\\
b_{n-1}=a_{0}+\ldots +a_{n-1}\in IM\cap J^{n}M=IJ^{n}M\mbox{ and
}a_{n-1}=b_{n-1}-b_{n-2}\\
a_{n}=-b_{n-1}\in IJ^{n}M\, .
\end{array}\right.$$

\noindent We can rewrite $z$ in $M[u,v]$ in the following manner:
\begin{eqnarray*}
&& z=a_{0}u^{n}+a_{1}u^{n-1}v+\ldots +a_{n-1}uv^{n-1}+a_{n}v^{n}=\\
&& =b_{0}u^{n}+(b_{1}-b_{0})u^{n-1}v+(b_{2}-b_{1})u^{n-2}v^{2}
+\ldots +\\ 
&&+(b_{n-2}-b_{n-3})u^{2}v^{n-2}+
(b_{n-1}-b_{n-2})uv^{n-1}+(-b_{n-1})v^{n}=\\
&& =\underbrace{(b_{0}u^{n-1}+b_{1}u^{n-2}v+b_{2}u^{n-3}v^{2}+\ldots
+b_{n-2}uv^{n-2}+b_{n-1}v^{n-1})}_{p(u,v)}(u-v):=p(u,v)(u-v)\, ,
\end{eqnarray*}
where $p(u,v)\in A[Iu,Jv]_{n-1}\cdot (IJM)=\mrees{I}{J}_{n-1}\cdot
(IJM)$. Since by hypothesis $IM\cap JM=IJM$, then ${\rm ker}(\sigma
_{1})=(IJM)(u-v)$, ${\rm ker}(\overline{\sigma}_{1})=0$ and
$z=p(u,v)(u-v)\in \mrees{I}{J}_{n-1}\cdot
(IJM)(u-v)=\mrees{I}{J}_{n-1}\cdot {\rm ker}\sigma _{1}$. Thus ${\rm
ker}\sigma _{n}=\mrees{I}{J}_{n-1}\cdot {\rm ker}\sigma _{1}$ for all
$n\geq 2$ and $s(E(\sigma ))=1$. By Lemma \ref{tecnic},
$s(E(\overline{\sigma}_{1}))\leq s(E(\sigma ))=1$. Therefore ${\rm
ker}(\overline{\sigma}_{n})=\magr{I}{J}_{n-1}\cdot {\rm
ker}(\overline{\sigma}_{1})=0$ for all $n\geq 2$ and
$\overline{\sigma}$ is an isomorphism. \fidemo

\medskip

\begin{proposicio}\label{normal} 
Let $A$ be a noetherian ring, $I$, $J$ two ideals of $A$ and $M$ a
finitely generated $A$-module. The following two conditions are
equivalent:
\begin{itemize}
\item[$(i)$] $\overline{\sigma}:\magrw{I}{J}{M}\rightarrow
\agrw{I+J}{M}$ is an isomorphism.
\item[$(ii)$] $I^{p}M\cap J^{q}M=I^{p}J^{q}M$ for all integers
$p,q\geq 1$.
\end{itemize}
\end{proposicio}

\noindent {\em Proof}. Remark that we can suppose $A$ is local. By
Lemma \ref{necesig}, $(ii)\Rightarrow (i)$. Let us see $(i)\Rightarrow
(ii)$, proving by double induction in $p,q\geq 1$ that
\begin{eqnarray*}
I^{p}M\cap J^{q}M\subset I^{p}(I+J)J^{q-1}M\cap (I+J)^{p}J^{q}M\, .
\end{eqnarray*}
Remark that if $I^{p}M\cap J^{q}M\subset I^{p}(I+J)J^{q-1}M$ for all
$p,q\geq 1$, then $I^{p}M\cap J^{q}M\subset I^{p+1}M+I^{p}J^{q}M$ and
$I^{p}M\cap J^{q}M\subset I^{p+1}M\cap
J^{q}M+I^{p}J^{q}M$. Recursively, and using $A$ is noetherian local
and $M$ is finitely generated, $I^{p}M\cap J^{q}M\subset (\cap _{r\geq
1}I^{p+r}M\cap J^{q}M)+I^{p}J^{q}M\subset (\cap _{n\geq 1}I^{n}M)+
I^{p}J^{q}M=I^{p}J^{q}M$, concluding $(ii)$. Take $q=1$. Let us prove
by induction in $p\geq 1$ that
\begin{eqnarray*}
I^{p}M\cap JM\subset I^{p}(I+J)M\cap (I+J)^{p}JM \, .
\end{eqnarray*}
For $p=1$, we apply Lemma \ref{necesig}, $(b)$, using the hypothesis
$\overline{\sigma}_{1}$ is an isomorphism. Suppose 
\begin{eqnarray*}
I^{p}M\cap JM\subset I^{p}(I+J)M\cap (I+J)^{p}JM 
\end{eqnarray*}
is true and let us prove
\begin{eqnarray*}
I^{p+1}M\cap JM\subset I^{p+1}(I+J)M\cap (I+J)^{p+1}JM \, .
\end{eqnarray*}
Then $I^{p+1}M\cap JM\subset I^{p}M\cap
JM\subset (I+J)^{p}JM$. Consider the short complex of $A$-modules:
\begin{eqnarray*}
I^{p+1}M\cap JM\buildrel \alpha \over\longrightarrow I^{p+1}M\oplus
(I+J)^{p}JM\buildrel \beta \over \longrightarrow (I+J)^{p+1}M\, ,
\end{eqnarray*} 
where $\alpha (a)=(a,-a)$ and $\beta (a,b)=a+b$. Remark that $\beta
\circ \alpha =0$, $\beta $ is surjective and that there exists a
natural epimorphism $\gamma $ of $A$-modules such that $\beta \circ
\gamma =\sigma _{p+1}$. If we tensor this short complex by $A/(I+J)$
we obtain:
\begin{eqnarray*}
&& (I^{p+1}M\cap JM)\otimes A/(I+J)\buildrel \overline{\alpha}
\over\longrightarrow I^{p+1}M/I^{p+1}(I+J)M\oplus
(I+J)^{p}JM/(I+J)^{p+1}JM \\ && I^{p+1}M/I^{p+1}(I+J)M\oplus
(I+J)^{p}JM/(I+J)^{p+1}JM \buildrel \overline{\beta} \over
\longrightarrow (I+J)^{p+1}M/(I+J)^{p+2}M\, ,
\end{eqnarray*} 
with $\overline{\beta}\circ \overline{\alpha}=0$. Since
$\overline{\sigma}_{p+1}=\overline{\beta}\circ \overline{\gamma}$ is
an isomorphism, then $\overline{\beta}$ is an isomorphism,
$\overline{\alpha}=0$ and
\begin{eqnarray*}
I^{p+1}M\cap JM\subset I^{p+1}(I+J)M\cap (I+J)^{p+1}JM\, .
\end{eqnarray*}
By the symmetry of the problem, the following inclusion is also true
for all $q\geq 1$:
\begin{eqnarray*}
IM\cap J^{q}M\subset (I+J)J^{q}M\cap I(I+J)^{q}M \, .
\end{eqnarray*}
In particular, if $I^{p}M\cap JM\subset I^{p}(I+J)M$ for all $p\geq
1$, then $I^{p}M\cap JM\subset I^{p+1}M+I^{p}JM$ and $I^{p}M\cap
JM\subset I^{p+1}M\cap JM+I^{p}JM$. Recursively, and using $A$ is
noetherian local and $M$ is finitely generated, $I^{p}M\cap JM\subset
(\cap _{r\geq 1}I^{p+r}M\cap JM)+I^{p}JM\subset (\cap _{n\geq
1}I^{n}M)+ I^{p}JM=I^{p}JM$ concluding $I^{p}M\cap JM=I^{p}JM$ for all
$p\geq 1$. Again, by the symmetry of the problem, $IM\cap
J^{q}M=IJ^{q}M$ for all $q\geq 1$. Now, suppose
\begin{eqnarray*}
I^{p}M\cap J^{q}M\subset I^{p}(I+J)J^{q-1}M\cap (I+J)^{p}J^{q}M 
\end{eqnarray*}
holds for all $p\geq 1$ and let us prove, by induction in $p\geq 1$,
that 
\begin{eqnarray*}
I^{p}M\cap J^{q+1}M\subset I^{p}(I+J)J^{q}M\cap (I+J)^{p}J^{q+1}M \, .
\end{eqnarray*}
Remark that if $I^{p}M\cap J^{q}M\subset I^{p}(I+J)J^{q-1}M$ for all
$p\geq 1$, then $I^{p}M\cap J^{q}M\subset I^{p+1}M+I^{p}J^{q}M$ and
$I^{p}M\cap J^{q}M\subset I^{p+1}M\cap
J^{q}M+I^{p}J^{q}M$. Recursively, and using $A$ is noetherian local
and $M$ is finitely generated, $I^{p}M\cap J^{q}M\subset (\cap _{r\geq
1}I^{p+r}M\cap J^{q}M)+I^{p}J^{q}M\subset (\cap _{n\geq 1}I^{n}M)+
I^{p}J^{q}M=I^{p}J^{q}M$ concluding $I^{p}M\cap J^{q}M=I^{p}J^{q}M$
for all $p\geq 1$. For $p=1$, we have to show:
\begin{eqnarray*}
IM\cap J^{q+1}M\subset I(I+J)J^{q}M\cap (I+J)J^{q+1}M \, .
\end{eqnarray*}
We have $IM\cap J^{q+1}M\subset IM\cap J^{q}M=IJ^{q}M$. Consider the
short complex of $A$-modules:
\begin{eqnarray*}
IM\cap J^{q+1}M\buildrel \alpha \over\longrightarrow I^{q+1}M\oplus
\ldots \oplus IJ^{q}M\oplus J^{q+1}M
\buildrel \sigma _{q+1} \over \longrightarrow (I+J)^{q+1}M\, ,
\end{eqnarray*} 
where $\alpha (a)=(0,\ldots ,0,a,-a)$. Remark that $\sigma _{q+1}
\circ \alpha =0$. If we tensor
this complex by $A/(I+J)$ we obtain $\overline{\sigma}_{q+1}\circ
\overline{\alpha}=0$. Since $\overline{\sigma}_{q+1}$ is an
isomorphism, then $\overline{\alpha}=0$ and 
\begin{eqnarray*}
IM\cap J^{q+1}M\subset
I(I+J)J^{q}M\cap (I+J)J^{q+1}M\, .
\end{eqnarray*}
Suppose now true
\begin{eqnarray*}
I^{p}M\cap J^{q+1}M\subset I^{p}(I+J)J^{q}M\cap (I+J)^{p}J^{q+1}M 
\end{eqnarray*}
and let us prove
\begin{eqnarray*}
I^{p+1}M\cap J^{q+1}M\subset I^{p+1}(I+J)J^{q}M\cap
(I+J)^{p+1}J^{q+1}M \, .
\end{eqnarray*}
Then $I^{p+1}M\cap J^{q+1}M\subset I^{p}M\cap J^{q+1}M\subset
(I+J)^{p}J^{q+1}M$ and $I^{p+1}M\cap J^{q+1}M\subset I^{p+1}M\cap
J^{q}M=I^{p+1}J^{q}M$. Consider the short complex of $A$-modules:
\begin{eqnarray*}
I^{p+1}M\cap J^{q+1}M\buildrel \alpha \over\longrightarrow
I^{p+q+1}M\oplus \ldots \oplus I^{p+1}J^{q}M\oplus
(I+J)^{p}J^{q+1}M\buildrel \beta \over \longrightarrow (I+J)^{p+q+1}M\, ,
\end{eqnarray*} 
where $\alpha (a)=(0,\ldots ,0,a,-a)$ and $\beta (a_{1},\ldots
,a_{q+2})=a_{1}+\ldots +a_{q+2}$. Remark that $\beta \circ \alpha =0$,
$\beta $ is surjective and that there exists a natural epimorphism
$\gamma $ of $A$-modules such that $\beta \circ \gamma =\sigma
_{p+q+1}$. If we tensor this complex by $A/(I+J)$ we obtain
$\overline{\beta} \circ \overline{\alpha}=0$. Since
$\overline{\sigma}_{p+q+1}=\overline{\beta}\circ \overline{\gamma}$ is
an isomorphism, then $\overline{\beta}$ is an isomorphism,
$\overline{\alpha}=0$ and
\begin{eqnarray*}
I^{p+1}M\cap J^{q+1}M\subset I^{p+1}(I+J)J^{q}M\cap
(I+J)^{p+1}J^{q+1}JM \, \, . \, \fidemo
\end{eqnarray*}

\begin{proposicio}\label{alaval}
Let $A$ be a commutative ring, $I$ an ideal of $A$ and $\lambda:
M\otimes N\rightarrow P$ an epimorphism of $A$-modules. Consider
$f:\reesw{I}{M}\otimes N\rightarrow \reesw{I}{P}$ and
$\overline{f}=f\otimes 1_{A/I}:\agrw{I}{M}\otimes N\rightarrow
\agrw{I}{P}$ the natural surjective graded morphisms of standard
modules. Then, for each integer $n\geq 2$, there exists an exact
sequence of $A$-modules $E(f)_{n+1}\rightarrow E(f)_{n}\rightarrow
E(\overline{f})_{n}\rightarrow 0$. In particular, if $A$ is
noetherian, $M,N,P$ are finitely generated and $\overline{f}$ is an
isomorphism, then $f$ is an isomorphism.
\end{proposicio}

\noindent {\em Proof}. For each integer $n\geq 1$, the natural
morphism $\tor{1}{}{A/I^{n}}{M}\otimes N\rightarrow
\tor{1}{}{A/I^{n}}{M\otimes N}$ and $\lambda :M\otimes N\rightarrow P$
define the following commutative diagram of exact rows:

\begin{picture}(330,85)(-5,0)

\put(100,60){\makebox(0,0){\mbox{\footnotesize
$\tor{1}{}{A/I^{n}}{M}\otimes N$}}}
\put(200,60){\makebox(0,0){\mbox{\footnotesize
$I^{n}\otimes M\otimes N$}}}
\put(300,60){\makebox(0,0){\mbox{\footnotesize
$I^{n}M\otimes N$}}}
\put(360,60){\makebox(0,0){$0$}}

\put(140,60){\vector(1,0){25}}
\put(235,60){\vector(1,0){30}}
\put(320,60){\vector(1,0){30}}

\put(40,20){\vector(1,0){30}}
\put(135,20){\vector(1,0){30}}
\put(235,20){\vector(1,0){30}}
\put(320,20){\vector(1,0){30}}

\put(30,20){\makebox(0,0){$0$}}
\put(100,20){\makebox(0,0){{\footnotesize 
$\tor{1}{}{A/I^{n}}{P}$}}}
\put(200,20){\makebox(0,0){{\footnotesize $I^{n}\otimes P$}}}
\put(300,20){\makebox(0,0){{\footnotesize $I^{n}P$}}}
\put(360,20){\makebox(0,0){$0$}}

\put(100,50){\vector(0,-1){20}}
\put(200,50){\vector(0,-1){20}}
\put(200,50){\vector(0,-1){16}}
\put(215,40){\makebox(0,0){\mbox{\footnotesize $1\otimes \lambda$}}} 
\end{picture}

\noindent We deduce an epimorphism $f_{n}:I^{n}M\otimes N\rightarrow
I^{n}P$. On the other hand, $\reesw{I}{M}\otimes M$ is a standard
$\rees{I}$-module and $f=\oplus _{n\geq 0}f_{n}:\reesw{I}{M}\otimes
N\rightarrow \reesw{I}{P}$ defines a surjective graded morphism of
standard $\rees{I}$-modules. If we tensor $f$ by $A/I$, we get
$\overline{f}:\agrw{I}{M}\otimes N\rightarrow \agrw{I}{P}$ a
surjective graded morphism of standard $\agr{I}$-modules.

Let $X$ be an $A$-module. The following is a commutative diagram of
exact columns with rows the last three nonzero terms of the complexes
$\mathcal{K}(\reesw{I}{X})_{n+1}$, $\mathcal{K}(\reesw{I}{X})_{n}$ and
$\mathcal{K}(\agrw{I}{X})_{n}$ (see Proposition 2.6 in \cite{planas2}
for more details):

\begin{picture}(330,125)(-25,0)

\put(0,100){\makebox(0,0){\mbox{\footnotesize
$\mathcal{K}(\reesw{I}{X})_{n+1}$}}}
\put(100,100){\makebox(0,0){\mbox{\footnotesize
${\bf \Lambda}_{2}I\otimes I^{n-1}X$}}}
\put(200,100){\makebox(0,0){\mbox{\footnotesize
$I\otimes I^{n}X$}}}
\put(300,100){\makebox(0,0){\mbox{\footnotesize
$I^{n+1}X$}}}
\put(370,100){\makebox(0,0){$0$}}

\put(142,100){\vector(1,0){16}}
\put(242,100){\vector(1,0){16}}
\put(330,100){\vector(1,0){30}}

\put(0,60){\makebox(0,0){\mbox{\footnotesize
$\mathcal{K}(\reesw{I}{X})_{n}$}}}
\put(100,60){\makebox(0,0){\mbox{\footnotesize
${\bf \Lambda}_{2}I\otimes I^{n-2}X$}}}
\put(200,60){\makebox(0,0){\mbox{\footnotesize
$I\otimes I^{n-1}X$}}}
\put(300,60){\makebox(0,0){\mbox{\footnotesize
$I^{n}X$}}}
\put(370,60){\makebox(0,0){$0$}}

\put(142,60){\vector(1,0){16}}
\put(242,60){\vector(1,0){16}}
\put(330,60){\vector(1,0){30}}

\put(148,20){\vector(1,0){10}}
\put(242,20){\vector(1,0){16}}
\put(330,20){\vector(1,0){30}}

\put(0,20){\makebox(0,0){\mbox{\footnotesize
$\mathcal{K}(\agrw{I}{X})_{n}$}}}
\put(100,20){\makebox(0,0){\mbox{\footnotesize 
${\bf \Lambda}_{2}I/I^{2}\otimes I^{n-2}X/I^{n-1}X$}}}
\put(200,20){\makebox(0,0){{\footnotesize 
$I/I^{2}\otimes I^{n-1}X/I^{n}X$}}}
\put(300,20){\makebox(0,0){\mbox{\footnotesize $I^{n}X/I^{n+1}X$}}}
\put(370,20){\makebox(0,0){$0$}}

\put(0,50){\vector(0,-1){20}}
\put(0,50){\vector(0,-1){16}}
\put(100,50){\vector(0,-1){20}}
\put(100,50){\vector(0,-1){16}}
\put(200,50){\vector(0,-1){20}}
\put(200,50){\vector(0,-1){16}}
\put(300,50){\vector(0,-1){20}}
\put(300,50){\vector(0,-1){16}}
\put(0,90){\vector(0,-1){20}}
\put(100,90){\vector(0,-1){20}}
\put(200,90){\vector(0,-1){20}}
\put(300,90){\vector(0,-1){20}}

\put(10,80){\makebox(0,0){\mbox{\footnotesize $u_{\cdot}$}}}
\put(110,80){\makebox(0,0){\mbox{\footnotesize $u_{2}$}}}
\put(210,80){\makebox(0,0){\mbox{\footnotesize $u_{1}$}}}
\put(310,80){\makebox(0,0){\mbox{\footnotesize $u_{0}$}}}
\put(10,40){\makebox(0,0){\mbox{\footnotesize $v_{\cdot}$}}}
\put(110,40){\makebox(0,0){\mbox{\footnotesize $v_{2}$}}}
\put(210,40){\makebox(0,0){\mbox{\footnotesize $v_{1}$}}}
\put(310,40){\makebox(0,0){\mbox{\footnotesize $v_{0}$}}}

\put(40,100){\makebox(0,0){\mbox{\footnotesize $\ldots$}}}
\put(40,60){\makebox(0,0){\mbox{\footnotesize $\ldots$}}}
\put(40,20){\makebox(0,0){\mbox{\footnotesize $\ldots$}}}
\put(150,110){\makebox(0,0){\mbox{\footnotesize $\partial _{2,n+1}$}}}
\put(150,70){\makebox(0,0){\mbox{\footnotesize $\partial _{2,n}$}}}
\put(150,10){\makebox(0,0){\mbox{\footnotesize $\partial _{2,n}$}}}
\put(250,110){\makebox(0,0){\mbox{\footnotesize $\partial _{1,n+1}$}}}
\put(250,70){\makebox(0,0){\mbox{\footnotesize $\partial _{1,n}$}}}
\put(250,10){\makebox(0,0){\mbox{\footnotesize $\partial _{1,n}$}}}
\end{picture}

\noindent In other words, $\mathcal{K}(\reesw{I}{X})_{n+1}\buildrel
u_{\cdot}\over \rightarrow \mathcal{K}(\reesw{I}{X})_{n}\buildrel
v_{\cdot}\over \rightarrow \mathcal{K}(\agrw{I}{X})_{n}\rightarrow 0$
is an exact sequence of complexes. It induces the morphisms in
homology: $H_{1}(\mathcal{K}(\reesw{I}{X})_{n+1})\buildrel u\over
\rightarrow H_{1}(\mathcal{K}(\reesw{I}{X})_{n})$ and
$H_{1}(\mathcal{K}(\reesw{I}{X})_{n})\buildrel v\over \rightarrow
H_{1}(\mathcal{K}(\agrw{I}{X})_{n})$. By Proposition 2.6 in
\cite{planas2}, $H_{1}(\mathcal{K}(\reesw{I}{X})_{n})=E(I;X)_{n}$ and
$H_{1}(\mathcal{K}(\agrw{I}{X})_{n})=E(\agrw{I}{X})_{n}$. Thus we have
$E(I;X)_{n+1}\buildrel u\over\rightarrow E(I;X)_{n}\buildrel v\over
\rightarrow E(\agrw{I}{X})_{n}$. Since $v_{\cdot}\circ u_{\cdot}=0$,
then $v\circ u=0$. Since $u_{0}$ is injective, then ${\rm ker}v\subset
{\rm im}u$. Since $H_{0}(\mathcal{K}(\reesw{I}{X})_{n+1})=0$, then $v$
is surjective. So $E(I;X)_{n+1}\buildrel u\over\rightarrow
E(I;X)_{n}\buildrel v\over \rightarrow E(\agrw{I}{X})_{n}\rightarrow
0$ is an exact sequence of $A$-modules. For $X=P$ we get the exact
sequence of $A$-modules: $E(I;P)_{n+1}\buildrel u\over\rightarrow
E(I;P)_{n}\buildrel v\over \rightarrow E(\agrw{I}{P})_{n}\rightarrow
0$. Take $X=M$ in $\mathcal{K}(\reesw{I}{X})_{n+1}\buildrel
u_{\cdot}\over \rightarrow \mathcal{K}(\reesw{I}{X})_{n}\buildrel
v_{\cdot}\over \rightarrow \mathcal{K}(\agrw{I}{X})_{n}\rightarrow 0$
and tensor it by $N$. Then we get the exact sequence of complexes
\begin{eqnarray*}
\mathcal{K}(\reesw{I}{M})_{n+1}\otimes N\buildrel
\alpha _{\cdot}=u_{\cdot}\otimes 1\over \longrightarrow 
\mathcal{K}(\reesw{I}{M})_{n}\otimes N\buildrel
\beta _{\cdot}=v_{\cdot}\otimes 1\over \longrightarrow 
\mathcal{K}(\agrw{I}{M})_{n}\otimes N\longrightarrow 0 \, .
\end{eqnarray*}
That is, we obtain the exact sequence:
\begin{eqnarray*}
\mathcal{K}(\reesw{I}{M}\otimes N)_{n+1}\buildrel
\alpha _{\cdot}\over \longrightarrow 
\mathcal{K}(\reesw{I}{M}\otimes N)_{n}\buildrel
\beta _{\cdot}\over \longrightarrow 
\mathcal{K}(\agrw{I}{M}\otimes N)_{n}\longrightarrow 0 \, ,
\end{eqnarray*}
which induces the morphisms in homology
\begin{eqnarray*}
H_{1}(\mathcal{K}(\reesw{I}{M}\otimes N)_{n+1})\buildrel \alpha \over
\rightarrow H_{1}(\mathcal{K}(\reesw{I}{M}\otimes N)_{n})
\buildrel \beta \over
\rightarrow H_{1}(\mathcal{K}(\agrw{I}{M}\otimes N)_{n})\, .
\end{eqnarray*}
Again, by Proposition 2.6 in \cite{planas2},
$H_{1}(\mathcal{K}(\reesw{I}{M}\otimes N)_{n})=E(\reesw{I}{M}\otimes
N)_{n}$ and $H_{1}(\mathcal{K}(\agrw{I}{M}\otimes
N)_{n})=E(\agr{I}\otimes M)_{n}$. Moreover, since $\beta _{\cdot}\circ
\alpha _{\cdot}=0$, then $\beta \circ \alpha =0$, and since
$H_{0}(\mathcal{K}(\reesw{I}{M}\otimes N)_{n+1})=0$, then $\beta $ is
an epimorphism. Thus we have
\begin{eqnarray*}
E(\reesw{I}{M}\otimes N)_{n+1}\buildrel \alpha \over \longrightarrow
E(\reesw{I}{M}\otimes N)_{n}\buildrel \beta \over \longrightarrow
E(\agrw{I}{M}\otimes N)_{n}\longrightarrow 0
\end{eqnarray*}
with $\beta \circ \alpha =0$ and $\beta$ surjective. Remark that since
we do not know if $\alpha _{0}=u_{0}\otimes 1$ is injective, we can
not deduce ${\rm ker}\beta \subset {\rm im}\alpha$.  On the other
hand, consider $g:{\bf S}(I)\otimes M\otimes N\rightarrow
\reesw{I}{M}\otimes N$ and $\overline{g}:{\bf S}(I/I^{2})\otimes
M\otimes N\rightarrow \agrw{I}{M}\otimes N$ the natural surjective
graded morphisms of standard modules, where ${\bf S}(I)$, ${\bf
S}(I/I^{2})$ stands for the symmetric algebras of $I$ and $I/I^{2}$,
respectively. By Lemma 2.3 in \cite{planas2}, for each $n\geq 2$, there
exists exact sequences of $A$-modules $E(g)_{n}\rightarrow E(f\circ
g)_{n}\rightarrow E(f)_{n}\rightarrow 0$ and
$E(\overline{g})_{n}\rightarrow E(\overline{f}\circ
\overline{g})_{n}\rightarrow E(\overline{f})_{n}\rightarrow 0$. In
other words, we have exact sequences
\begin{eqnarray*}
&& E(\reesw{I}{M}\otimes N)_{n}\rightarrow
E(\reesw{I}{P})_{n}\rightarrow E(f)_{n}\rightarrow 0\mbox{ and } \\&&
E(\agrw{I}{M}\otimes N)_{n}\rightarrow E(\agrw{I}{P})_{n}\rightarrow
E(\overline{f})\rightarrow 0\, .
\end{eqnarray*}
Consider the following commutative diagram of exact columns:

\begin{picture}(330,125)(0,0)

\put(100,100){\makebox(0,0){\mbox{\footnotesize
$E(\reesw{I}{M}\otimes N)_{n+1}$}}}
\put(200,100){\makebox(0,0){\mbox{\footnotesize
$E(\reesw{I}{M}\otimes N)_{n}$}}}
\put(300,100){\makebox(0,0){\mbox{\footnotesize
$E(\agrw{I}{M}\otimes N)_{n}$}}}
\put(370,100){\makebox(0,0){$0$}}

\put(142,100){\vector(1,0){16}}
\put(242,100){\vector(1,0){16}}
\put(335,100){\vector(1,0){25}}

\put(100,60){\makebox(0,0){\mbox{\footnotesize
$E(\reesw{I}{P})_{n+1}$}}}
\put(200,60){\makebox(0,0){\mbox{\footnotesize
$E(\reesw{I}{P})_{n}$}}}
\put(300,60){\makebox(0,0){\mbox{\footnotesize
$E(\agrw{I}{P})_{n}$}}}
\put(370,60){\makebox(0,0){$0$}}

\put(142,60){\vector(1,0){16}}
\put(242,60){\vector(1,0){16}}
\put(330,60){\vector(1,0){30}}

\put(100,20){\makebox(0,0){\mbox{\footnotesize 
$E(f)_{n+1}$}}}
\put(200,20){\makebox(0,0){{\footnotesize 
$E(f)_{n}$}}}
\put(300,20){\makebox(0,0){\mbox{\footnotesize 
$E(\overline{f})_{n}$}}}

\put(100,50){\vector(0,-1){20}}
\put(100,50){\vector(0,-1){16}}
\put(200,50){\vector(0,-1){20}}
\put(200,50){\vector(0,-1){16}}
\put(300,50){\vector(0,-1){20}}
\put(300,50){\vector(0,-1){16}}
\put(100,90){\vector(0,-1){20}}
\put(200,90){\vector(0,-1){20}}
\put(300,90){\vector(0,-1){20}}

\put(150,110){\makebox(0,0){\mbox{\footnotesize $\alpha$}}}
\put(150,70){\makebox(0,0){\mbox{\footnotesize $u$}}}
\put(250,110){\makebox(0,0){\mbox{\footnotesize $\beta$}}}
\put(250,70){\makebox(0,0){\mbox{\footnotesize $v$}}}
\end{picture}

\noindent The commutativity induces two morphisms $\xi
:E(f)_{n+1}\rightarrow E(f)_{n}$ and $\mu :E(f)_{n}\rightarrow
E(\overline{f})_{n}$. Since $v\circ u=0$, then $\mu \circ \xi
=0$. Since $v$ is surjective, then $\mu$ is surjective too. Since
$\beta$ is surjective and the middle row is exact, then ${\rm
ker}\mu\subset {\rm im}\xi $. Therefore,
\begin{eqnarray*}
E(f)_{n+1}\buildrel \xi \over \longrightarrow E(f)_{n}\buildrel
\mu \over \longrightarrow E(\overline{f})_{n}\longrightarrow 0
\end{eqnarray*}
is an exact sequence of $A$ modules. Finally, if $A$ is noetherian and
$M,N$ and $P$ are finitely generated, then $E(f)_{n}=0$ for $n\gg 0$
big enough. \fidemo

\medskip

\begin{teorema}\label{trans}
Let $A$ be a noetherian ring, $I$, $J$ two ideals of $A$ and $M$ a
finitely generated $A$-module. The following two conditions are
equivalent:
\begin{itemize}
\item[$(i)$] $\overline{\varphi}:\agr{I}\otimes
\agrw{J}{M}\rightarrow \agrw{I+J}{M}$ is an isomorphism.
\item[$(ii)$] $\tor{1}{}{A/I^{p}}{\reesw{J}{M}}=0$ and
$\tor{1}{}{A/I^{p}}{\agrw{J}{M}}=0$ for all integers $p\geq 1$.
\end{itemize}
In particular, $\agr{I}\otimes \agr{J}\simeq \agr{I+J}$ if and only if
$\tor{1}{}{A/I^{p}}{A/J^{q}}=0$ and $\tor{2}{}{A/I^{p}}{A/J^{q}}=0$
for all integers $p,q\geq 1$.
\end{teorema}

\noindent {\em Proof}. Remark that $\tor{1}{}{A/I^{p}}{J^{q}M}={\rm
ker}(\pi _{p,q}:I^{p}\otimes J^{q}M\rightarrow
I^{p}J^{q}M)$. Moreover, under the hypothesis
$\tor{1}{}{A/I^{p}}{\reesw{J}{M}}=0$ for all $p\geq 1$, then the
following two conditions are equivalent:
\begin{itemize}
\item $\tor{1}{}{A/I^{p}}{\agrw{J}{M}}=0$ for all $p\geq 1$. 
\item $I^{p}M\cap J^{q}M=I^{p}J^{q}M$ for all $p,q\geq 1$.
\end{itemize}
Suppose $(ii)$ holds, i.e., $\tor{1}{}{A/I^{p}}{J^{q}M}=0$ and
$I^{p}M\cap J^{q}M=I^{p}J^{q}M$ for all $p,q\geq 1$. Then, $\pi
:\rees{I}\otimes \reesw{J}{M}\rightarrow \mreesw{I}{J}{M}$ is an
isomorphism and, by Lemma \ref{necesig},
$\overline{\sigma}:\magrw{I}{J}{M}\rightarrow \agrw{I+J}{M} $ is an
isomorphism. Thus $\overline{\varphi}=\overline{\sigma}\circ
\overline{\pi}$ is an isomorphism and $(i)$ holds. Let us now prove
$(i)\Rightarrow (ii)$. If $\overline{\varphi}=\overline{\sigma}\circ
\overline{\pi}$ is an isomorphism , then $\overline{\sigma}$ and
$\overline{\pi}$ are two isomorphisms. By Proposition \ref{normal},
$\overline{\sigma}$ an isomorphism implies $I^{p}M\cap
J^{q}M=I^{p}J^{q}M$ for all $p,q\geq 1$. In particular,
\begin{eqnarray*}
&&\reesw{I}{J^{q}M/J^{q+1}M}_{p}
=\frac{I^{p}J^{q}M+J^{q+1}M}{J^{q+1}M}=
\frac{I^{p}J^{q}M}{I^{p}J^{q}M\cap
J^{q+1}M}=\frac{I^{p}J^{q}M}{I^{p}J^{q+1}M}=
\agrw{J}{I^{p}M}_{q}\mbox{ and}\\
&&\agrw{I}{J^{q}M/J^{q+1}M}_{p}=
\frac{I^{p}J^{q}M+J^{q+1}M}{I^{p+1}J^{q}M+J^{q+1}M}=
\frac{I^{p}J^{q}M}{(I+J)I^{p}J^{q}M}=\magrw{I}{J}{M}_{p,q}\, .
\end{eqnarray*}
Fix $q\geq 1$. Since $\overline{\pi} _{p,q}:\agr{I}_{p}\otimes
\agrw{J}{M}_{q}\rightarrow \magrw{I}{J}{M}_{p,q}$ is an isomorphism
for all $p\geq 1$ and
$\magrw{I}{J}{M}_{p,q}=\agrw{I}{J^{q}M/J^{q+1}M}_{p}$, then
$\overline{\pi}_{*,q}:\agr{I}\otimes J^{q}M/J^{q+1}M\rightarrow
\agrw{I}{J^{q}M/J^{q+1}M}$ is an isomorphism for all $q\geq 1$. By
Proposition \ref{alaval}, we have $\rees{I}\otimes
J^{q}M/J^{q+1}M\rightarrow \reesw{I}{J^{q}M/J^{q+1}M}$ is an
isomorphism for all $q\geq 1$. In other words, $I^{p}\otimes
\agrw{J}{M}\rightarrow \agrw{J}{I^{p}M}$ is an isomorphism for all
$p\geq 1$ (since
$\reesw{I}{J^{q}M/J^{q+1}M}_{p}=\agrw{J}{I^{p}M}_{q}$). By
Proposition \ref{alaval}, $I^{p}\otimes \reesw{J}{M}\rightarrow
\reesw{J}{I^{p}M}$ is an isomorphism for all $p\geq 1$. So $\pi
:\rees{I}\otimes \reesw{J}{M}\rightarrow \mreesw{I}{J}{M}$ is an
isomorphism and $\tor{1}{}{A/I^{p}}{\mrees{J}{M}}=0$ for all $p\geq
1$. \fidemo

\section{Some examples}

\begin{exemple}{\rm Let $A$ be a noetherian local ring, $I$, $J$ two
ideals of $A$ and $M$ a finitely generated $A$-module. If $I=(x)$ is
principal and $x$ $A$-regular, then $\overline{\varphi}:\agr{I}\otimes
\agrw{J}{M}\rightarrow \agrw{I+J}{M}$ is an isomorphism if and only if
$x$ is a nonzero divisor in $\reesw{J}{M}$ and in
$\agrw{J}{M}$. Indeed, let $\mathcal{K}(y;N)$ denote the Koszul
complex of a sequence of elements $y=y_{1},\ldots ,y_{m}$ of $A$ with
respect to an $A$-module $N$ and let $H_{i}(y;N)$ denote its $i$-th
Koszul homology group. Then
$\tor{1}{}{A/I}{N}=H_{1}(\mathcal{K}(x;A)\otimes N)=H_{1}(x;M) =0$ if
and only if $x$ is a non-zerodivisor in $N$.}\end{exemple}

\begin{exemple}{\rm Let $A$ be a noetherian local ring and let $I=(x)$
and $J=(y)$ be two principal ideals of $A$. If $(0:x)\subset (y)$ and
$(0:y)\subset (x)$, then $\overline{\varphi}:\agr{I}\otimes
\agr{J}\rightarrow \agr{I+J}$ is an isomorphism if and only if $x,y$
is an $A$-regular sequence.
}\end{exemple}

\begin{exemple}{\rm Let $R$ be a noetherian local ring and let $z,t$
be an $R$-regular sequence. Let $A=R/(zt)$, $x=z+(zt)$, $y=t+(zt)$,
$I=(x)$ and $J=(y)$. Then $\overline{\sigma}:\magr{I}{J}\rightarrow
\agr{I+J}$ is an isomorphism, but $\overline{\pi}:\agr{I}\otimes
\agr{J}\rightarrow \magr{I}{J}$ is not an isomorphism.  }\end{exemple}

An example of a pair of ideals $I$, $J$ with the property
$\tor{1}{}{A/I^{p}}{A/J^{q}}$ for all integers $p,q\geq 1$ arises from
a product of affine varietes (see \cite{vasconcelos}, pages 130 to
136, and specially Proposition 5.5.7). The next result is well known
(see, for instance, \cite{hio}). We give here a proof for the sake of
completeness.

\begin{proposicio}\label{flatness}
Let $A$ be a noetherian local ring, $I$ and $J$ two ideals of $A$ and
$M$ a finitely generated $A$-module. Let $x=x_{1},\ldots ,x_{r}$ be a
system of generators of $I$ and $y=y_{1},\ldots ,y_{r}$,
$y_{i}=\overline{x}_{i}=x_{i}+J$, a system of generators of the ideal
$\overline{I}=I+J/J$ of the quotient ring $\overline{A}=A/J$. If
$\agr{J}$ and $\agrw{J}{M}$ are free $\overline{A}$-modules and $y$ is
an $\overline{A}$-regular sequence in $\overline{I}$, then $x$ is an
$A$-regular sequence in $I$ and then
$\overline{\varphi}:\agr{I}\otimes \agrw{J}{M}\rightarrow
\agrw{I+J}{M}$ is an isomorphism.
\end{proposicio}

\noindent {\em Proof}. Since, for all $q\geq 1$, $J^{q}M/J^{q+1}M$ is
$\overline{A}$-free and $y$ is an $\overline{A}$-regular sequence,
then
\begin{eqnarray*}
0=\tor{1}{\overline{A}}{\overline{A}/\overline{I}}{J^{q}M/J^{q+1}M}
=H_{1}(\mathcal{K}(y;\overline{A})\otimes J^{q}M/J^{q+1}M)
=H_{1}(y;J^{q}M/J^{q+1}M)\, .
\end{eqnarray*}
So $y$ is a $J^{q}M/J^{q+1}M$-regular sequence in $\overline{I}$
for all $q\geq 1$. In particular, $x$ is a $J^{q}M/J^{q+1}M$-regular
sequence in $I$ and $H_{1}(x;J^{q}M/J^{q+1}M)=0$ for all $q\geq
1$. Using the long exact sequences in homology associated to the short
exact sequences of $A$-modules $0\rightarrow
J^{q}M/J^{q+1}M\rightarrow M/J^{q+1}M\rightarrow M/J^{q}M\rightarrow
0$, we deduce $H_{1}(x;M/J^{q}M)=0$ and $x$ is an $M/J^{q}M$-regular
sequence in $I$ for all $q\geq 1$. In particular, $x$ is an
$M$-regular sequence in $I$. Analogously, but using the hypothesis
$\agr{J}$ is $\overline{A}$-free, we deduce $x$ is an $A$-regular
sequence in $I$. Therefore
\begin{eqnarray*}
&& \tor{i}{}{A/I}{M}=H_{i}(\mathcal{K}(x;A)\otimes M)
=H_{i}(\mathcal{K}(x;M))=0\mbox{ and}\\ &&
\tor{i}{}{A/I}{M/J^{q}M}=H_{i}(\mathcal{K}(x;A)\otimes M/J^{q}M)
=H_{i}(\mathcal{K}(x;M/J^{q}M))=0\, .
\end{eqnarray*}
Using the long exact sequences in homology associated to the short
exact sequences
\begin{eqnarray*}
0\rightarrow J^{q}M\rightarrow M\rightarrow
M/J^{q}M\rightarrow 0 \mbox{ and }
0\rightarrow J^{q}M/J^{q+1}M\rightarrow
M/J^{q+1}M\rightarrow M/J^{q}M\rightarrow 0\, ,
\end{eqnarray*} we deduce
$\tor{1}{}{A/I}{\reesw{J}{M}}=0$ and $\tor{1}{}{A/I}{\agrw{J}{M}}=0$.
Since $I^{p}/I^{p+1}$ is $A/I$-free, then
$\tor{1}{}{I^{p}/I^{p+1}}{\reesw{J}{M}}
=\tor{1}{}{A/I}{\reesw{J}{M}}\otimes I^{p}/I^{p+1}=0$ and
$\tor{1}{}{I^{p}/I^{p+1}}{\agrw{J}{M}}
=\tor{1}{}{A/I}{\agrw{J}{M}}\otimes I^{p}/I^{p+1}=0$. Applying the
long exact sequences in homology to the short exact sequences
$0\rightarrow I^{p}/I^{p+1}\rightarrow A/I^{p+1}\rightarrow
A/I^{p}\rightarrow 0$, we deduce $\tor{1}{}{A/I^{p}}{\reesw{J}{M}}=0$
and $\tor{1}{}{A/I^{p}}{\agrw{J}{M}}=0$ for all $p\geq 1$. \fidemo

\section{Relation type of tensor products}

\begin{lema}\label{preli}
Let $U$ be a standard $A$-algebra and $F$ a standard $U$-module. If
$M$ is an $A$-module, then $F\otimes M$ is a standard $U$-module and
${\rm rt}(F\otimes M)\leq {\rm rt}(F)$. If $\lambda :M\rightarrow N$
is an epimorphism of $A$-modules, then $1\otimes \lambda :F\otimes
M\rightarrow F\otimes N$ is a surjective graded morphism of standard
$U$-modules. Moreover, for each integer $n\geq 1$, ${\rm
ker}(1_{F_{n}}\otimes \lambda )=U_{1}\cdot {\rm
ker}(1_{F_{n-1}}\otimes \lambda )$. In particular, for each $n\geq 1$,
there exists an epimorphism of $A$-modules $E(F\otimes
M)_{n}\rightarrow E(F\otimes N)_{n}$ and ${\rm rt}(F\otimes N)\leq
{\rm rt}(F\otimes M)$.
\end{lema}

\noindent {\em Proof}. Clearly $F\otimes M$ is a standard $U$-module
and $1\otimes \lambda :F\otimes M\rightarrow F\otimes N$ is a
surjective graded morphism of standard $U$-modules. By Proposition 2.6
in \cite{planas2}, for each $n\geq {\rm rt}(F)+1$, the following
sequence is exact:
\begin{eqnarray*}
{\bf \Lambda}_{2}(U_{1})\otimes F_{n-2}\rightarrow U_{1}\otimes
F_{n-1}\rightarrow F_{n}\rightarrow 0\, .
\end{eqnarray*}
If we tensor it by $M$, we obtain the exact sequence
\begin{eqnarray*}
{\bf \Lambda}_{2}(U_{1})\otimes F_{n-2}\otimes M\rightarrow U_{1}\otimes
F_{n-1}\otimes M\rightarrow F_{n}\otimes M\rightarrow 0\, ,
\end{eqnarray*}
for all $n\geq {\rm rt}(F)+1$. Thus $E(F\otimes M)_{n}=0$ for all
$n\geq {\rm rt}(F)+1$ and ${\rm rt}(F\otimes M)\leq {\rm
rt}(F)$. Consider the following commutative diagram of exact columns
and rows:

\begin{picture}(330,125)(0,0)

\put(200,100){\makebox(0,0){\mbox{\footnotesize
$({\rm ker}\partial _{n})\otimes M$}}}
\put(300,100){\makebox(0,0){\mbox{\footnotesize
$({\rm ker}\partial _{n})\otimes N$}}}
\put(370,100){\makebox(0,0){$0$}}

\put(242,100){\vector(1,0){16}}
\put(330,100){\vector(1,0){30}}

\put(100,60){\makebox(0,0){\mbox{\footnotesize 
$U_{1}\otimes {\rm ker}(1_{F_{n-1}}\otimes \lambda )$}}}
\put(200,60){\makebox(0,0){\mbox{\footnotesize 
$U_{1}\otimes F_{n-1}\otimes M$}}}
\put(300,60){\makebox(0,0){\mbox{\footnotesize 
$U_{1}\otimes F_{n-1}\otimes N$}}}
\put(370,60){\makebox(0,0){$0$}}

\put(142,60){\vector(1,0){16}}
\put(242,60){\vector(1,0){16}}
\put(333,60){\vector(1,0){27}}
\put(40,20){\vector(1,0){25}}
\put(138,20){\vector(1,0){24}}
\put(238,20){\vector(1,0){24}}
\put(330,20){\vector(1,0){30}}

\put(30,20){\makebox(0,0){$0$}}
\put(100,20){\makebox(0,0){\mbox{\footnotesize 
${\rm ker}(1_{F_{n}}\otimes \lambda )$}}}
\put(200,20){\makebox(0,0){{\footnotesize 
$F_{n}\otimes M$}}}
\put(300,20){\makebox(0,0){\mbox{\footnotesize 
$F_{n}\otimes N$}}}
\put(370,20){\makebox(0,0){$0$}}

\put(100,50){\vector(0,-1){20}}
\put(200,50){\vector(0,-1){20}}
\put(200,50){\vector(0,-1){16}}
\put(300,50){\vector(0,-1){20}}
\put(300,50){\vector(0,-1){16}}
\put(200,90){\vector(0,-1){20}}
\put(300,90){\vector(0,-1){20}}

\put(250,110){\makebox(0,0){\mbox{\footnotesize $1\otimes \lambda $}}}
\put(250,70){\makebox(0,0){\mbox{\footnotesize $1\otimes 1\otimes
\lambda $}}}
\put(222,40){\makebox(0,0){\mbox{\footnotesize $\partial _{n}\otimes
1_{M}$}}} 
\put(322,40){\makebox(0,0){\mbox{\footnotesize $\partial _{n}\otimes
1_{N}$}}} 
\end{picture}

\noindent Using a diagram chasing argument, one deduces ${\rm
ker}(1_{F_{n}}\otimes \lambda)=U_{1}\cdot {\rm ker}(1_{F_{n-1}}\otimes
\lambda )$ for all $n\geq 1$. If $g:X\rightarrow F\otimes M$ is a
symmetric presentation of $F\otimes M$, then, by Lemma 2.3 in
\cite{planas2}, there exists an exact sequence of $A$-modules
$E(g)_{n}\rightarrow E((1\otimes \lambda)\circ g)_{n}\rightarrow
E(1\otimes \lambda )_{n}\rightarrow 0$ for all $n\geq 1$. But
$E(g)_{n}=E(F\otimes M)_{n}$, $E((1\otimes \lambda )\circ
g)_{n}=E(F\otimes N)_{n}$ and $E(1\otimes \lambda )_{n}=0$ for all
$n\geq 1$. Thus $E(F\otimes M)_{n}\rightarrow E(F\otimes N)_{n}$ is
surjective for all $n\geq 1$ and ${\rm rt}(F\otimes N)\leq {\rm
rt}(F\otimes M)$. \fidemo

\begin{teorema}\label{rt}
Let $A$ be a commutative ring, $U$ and $V$ two standard $A$-algebras
and $F$ a standard $U$-module and $G$ a standard $V$-module. Then
$U\otimes V$ is a standard $A$-algebra, $F\otimes G$ is a standard
$U\otimes V$-module and ${\rm rt}(F\otimes G)\leq {\rm
max}({\rm rt}(F),{\rm rt}(G))$.
\end{teorema}

\noindent {\em Proof}. Clearly $U\otimes V$ is a standard $A$-algebra
and $F\otimes G$ is a standard $U\otimes V$-module. Take $\varphi
:X\rightarrow F$ and $\psi :Y\rightarrow G$ two symmetric
presentations of $F$ and $G$, respectively. Then $\varphi \otimes \psi
:X\otimes Y\rightarrow F\otimes G$ is a symmetric presentation of
$F\otimes G$. Since $\varphi \otimes \psi =(\varphi \otimes
1_{G})\circ (1_{X}\otimes \psi )$, then, for each integer $n\geq 2$,
there exists an exact sequence of $A$-modules
\begin{eqnarray*}
E(1_{X}\otimes \psi)_{n}\rightarrow E(\varphi \otimes \psi
)_{n}\rightarrow E(\varphi \otimes 1_{G})_{n}\rightarrow 0\, .
\end{eqnarray*}
Since $\psi :Y\rightarrow G$ is a symmetric presentation of $G$, then
$1_{X_{0}}\otimes \psi :X_{0}\otimes Y\rightarrow X_{0}\otimes G$ is a
symmetric presentation of $X_{0}\otimes G$ and $E(X_{0}\otimes G)_{n}=
E(1_{X_{0}}\otimes \psi )_{n}$. Using Lemma \ref{preli}, ${\rm
ker}(1_{X_{i}}\otimes \psi _{n-i})=U_{1}\cdot {\rm
ker}(1_{X_{i-1}}\otimes \psi _{n-i})$ for all $i\geq 1$. Then
\begin{eqnarray*}
&& E(1_{X}\otimes \psi )_{n}=\frac{{\rm ker}(1_{X}\otimes \psi
)_{n}}{(U\otimes V)_{1}\cdot {\rm ker}(1_{X}\otimes \psi )_{n-1}}=\\
&& \frac{\oplus _{i=0}^{n}{\rm ker}(1_{X_{i}}\otimes \psi
_{n-i})}{\left( \oplus _{i=0}^{n-1}U_{1}\cdot {\rm
ker}(1_{X_{i}}\otimes \psi _{n-i})\right) + \left( \oplus
_{i=0}^{n-1}V_{1}\cdot {\rm ker}(1_{X_{i}}\otimes \psi _{n-i})\right)
} = \\
&&\frac{{\rm ker}(1_{X_{0}}\otimes \psi _{n})}{V_{1}\cdot {\rm
ker}(1_{X_{0}}\otimes \psi _{n-1})}\oplus \frac{{\rm
ker}(1_{X_{1}}\otimes \psi _{n-1})}{ U_{1}\cdot {\rm
ker}(1_{X_{0}}\otimes \psi _{n-1})+ V_{1}\cdot {\rm
ker}(1_{X_{1}}\otimes \psi _{n-2})}\oplus \ldots \oplus \\
&& \frac{{\rm ker}(1_{X_{n-1}}\otimes \psi _{1})}{U_{1}\cdot {\rm
ker}(1_{X_{n-2}}\otimes \psi _{1})+ V_{1}\cdot {\rm
ker}(1_{X_{n-1}}\otimes \psi _{0})}\oplus \frac{{\rm
ker}(1_{X_{n}}\otimes \psi _{0})}{U_{1}\cdot {\rm
ker}(1_{X_{n-1}}\otimes \psi _{0})} = E(1_{X_{0}}\otimes \psi )_{n}\, .
\end{eqnarray*}
Therefore $E(1_{X}\otimes \psi )_{n}=E(1_{X_{0}}\otimes \psi
)_{n}=E(X_{0}\otimes G)_{n}$ for all $n\geq 1$.  Analogously,
$E(\varphi \otimes 1_{G})_{n}=E(\varphi \otimes
1_{G_{0}})_{n}=E(F\otimes G_{0})_{n}$ for all $n\geq 1$. Hence there
exists an exact sequence of $A$-modules
\begin{eqnarray*}
E(X_{0}\otimes G)_{n}\rightarrow E(F\otimes G)_{n}\rightarrow
E(F\otimes G_{0})_{n}\rightarrow 0
\end{eqnarray*}
for all $n\geq 2$ and, by Lemma \ref{preli}, ${\rm rt}(F\otimes G)\leq
{\rm max}({\rm rt}(F\otimes G_{0}), {\rm rt}(X_{0}\otimes G))\leq {\rm
max}({\rm rt}(F),{\rm rt}(G))$. \fidemo

\begin{observacio}{\rm 
Let $A$ be a commutative ring and let $U$ and $V$ be two standard
$A$-algebras. If $\tor{1}{A}{U}{V}=0$, then $E(U\otimes V)=E(U)\oplus
E(V)$. This follows from the characterization $E(U)=H_{1}(A,U,A)$ (see
Remark~2.3 in \cite{planas1}) and Proposition 19.3 in \cite{andre}.
}\end{observacio}

\section{Uniform bounds}

\begin{lema}\label{casli}
Let $(A,\mathfrak{m})$ be a noetherian local ring and $M$ be a
finitely generated $A$-module. Let $\mathfrak{p}$ a prime ideal of $A$
such that $A/\mathfrak{p}$ is regular local and $\agr{\mathfrak{p}}$
and $\agrw{\mathfrak{p}}{M}$ are free $A/\mathfrak{p}$-modules. Then
${\rm rt}(\mathfrak{m};M)\leq {\rm rt}(\mathfrak{p};M)$.
\end{lema}

\noindent {\em Proof}. Since $A/\mathfrak{p}$ is regular local, there
exists a sequence of elements $x=x_{1},\ldots ,x_{r}$ in $A$ such that
$y=y_{1},\ldots ,y_{r}$, defined by $y_{i}=x_{i}+\mathfrak{p}$, is a
system of generators of $\mathfrak{m}/\mathfrak{p}$ and an
$\overline{A}$-regular sequence. Let $I$ be the ideal of $A$ generated
by $x$. In particular,
$I+\mathfrak{p}/\mathfrak{p}=\mathfrak{m}/\mathfrak{p}$ and
$I+\mathfrak{p}=\mathfrak{m}$. By Proposition~\ref{flatness}, $x$ is
an $A$-regular sequence and
$\tor{1}{}{A/I^{p}}{\reesw{\mathfrak{p}}{M}}=0$ and
$\tor{1}{}{A/I^{p}}{\agrw{\mathfrak{p}}{M}}=0$ for all $p\geq 1$.  By
Theorem \ref{trans}, $\overline{\varphi}:\agr{I}\otimes
\agrw{\mathfrak{p}}{M}\rightarrow \agrw{\mathfrak{m}}{M}$ is an
isomorphism.  By Theorem \ref{rt}, ${\rm
rt}(\agrw{\mathfrak{m}}{M})\leq {\rm max}({\rm rt}(\agr{I}),{\rm
rt}(\agrw{\mathfrak{p}}{M}))$. By Remark 2.7 in \cite{planas2}, ${\rm
rt}(\agrw{J}{M})={\rm rt}(J;M)$ for any ideal $J$ of $A$. Since $I$ is
generated by a regular sequence, then ${\rm rt}(I)=1$ (see, for
instance, \cite{vasconcelos} page 30). Thus ${\rm
rt}(\mathfrak{m};M)\leq {\rm rt}(\mathfrak{p};M)$. \fidemo

\medskip

Next result is a slight generalization of a well known Theorem of
Duncan and O'Carroll \cite{do}. In fact the proof of our theorem is
directly inspired in their. We sketch it here for the sake of
completeness.

\begin{teorema}\label{boundmax}
Let $A$ be an excellent (or $J-2$) ring and let $M$ be a finitely
generated $A$-module.  Then there exists an integer $s\geq 1$ such
that, for all maximal ideals $\mathfrak{m}$ of $A$, the relation type
of $\mathfrak{m}$ with respect to $M$ satisfies ${\rm
rt}(\mathfrak{m};M)\leq s$.
\end{teorema}

\noindent {\em Proof}. For every $\mathfrak{p}\in {\rm Spec}(A)$, let
us construct a non-empty open subset $U(\mathfrak{p})$ of
$V(\mathfrak{p})= \{ \mathfrak{q}\in {\rm Spec}(A)\mid
\mathfrak{q}\supseteq \mathfrak{p}\} \simeq {\rm
Spec}(A/\mathfrak{p})$.  Remark that $A/\mathfrak{p}$ is a noetherian
domain, $\agr{\mathfrak{p}}$ is a finitely generated
$A/\mathfrak{p}$-algebra and $\agrw{\mathfrak{p}}{M}$ is a finitely
generated $\agr{\mathfrak{p}}$-module. By Generic Flatness (Theorem
22.A in \cite{matsumura}), there exist $f,g\in A-\mathfrak{p}$ such
that $\agr{\mathfrak{p}}_{f}$ is an $(A/\mathfrak{p})_{f}$-free module
and $\agrw{\mathfrak{p}}{M}_{g}$ is an $(A/\mathfrak{p})_{g}$-free
module. Since $A$ is $J-2$, the set ${\rm
Reg}(A/\mathfrak{p})=\{\mathfrak{q}\in V(\mathfrak{p})\mid
(A/\mathfrak{p})_{\mathfrak{q}} \mbox{ is regular local}\} $ is a
non-empty open subset of $V(\mathfrak{p})$. Define $U(\mathfrak{p})$
as the intersection $D(f)\cap D(g)\cap {\rm Reg}(A/\mathfrak{p})=\{
\mathfrak{q}\in V(\mathfrak{p})\mid \mathfrak{q}\not\ni f \, ,
\mathfrak{q}\not\ni g\, , (A/\mathfrak{p})_{\mathfrak{q}} \mbox{ is
regular local}\} $, which is a non-empty open subset of
$V(\mathfrak{p})$. Remark that for all $\mathfrak{q}\in
U(\mathfrak{p})$, $(A/\mathfrak{p})_{\mathfrak{q}}$ is regular local
and $\agr{\mathfrak{p}}_{\mathfrak{q}}$ and
$\agrw{\mathfrak{p}}{M}_{\mathfrak{q}}$ are free
$\agr{\mathfrak{p}}_{\mathfrak{q}}$-modules. By Lemma \ref{casli},
${\rm rt}(\mathfrak{q}A_{\mathfrak{q}};M_{\mathfrak{q}})\leq {\rm
rt}(\mathfrak{p}A_{\mathfrak{q}};M_{\mathfrak{q}})\leq {\rm
rt}(\mathfrak{p};M)$ for all $\mathfrak{q}\in U(\mathfrak{p})$. In
particular, ${\rm rt}(\mathfrak{m};M)\leq {\rm rt}(\mathfrak{p};M)$
for all maximal ideals $\mathfrak{m}\in U(\mathfrak{p})$. For each
minimal prime $\mathfrak{p}_{i}$ of $A$, let
$V(\mathfrak{p}_{i})-U(\mathfrak{p}_{i})=V(\mathfrak{p}_{i,1})\cup
\ldots \cup V(\mathfrak{p}_{i,r_{i}})$ be the decomposition into
irreducible closed subsets of the proper closed subset
$V(\mathfrak{p}_{i})-U(\mathfrak{p}_{i})$, $\mathfrak{p}_{i,j}\in {\rm
Spec}(A)$, $\mathfrak{p}_{i,j}\varsupsetneq \mathfrak{p}_{i}$.  Since
$A$ is noetherian, ${\rm Spec}(A)$ can be covered by finitely many
locally closed sets of type $U(\mathfrak{p})$, i.e., there exists a
finite number of prime ideals $\mathfrak{q}_{1},\ldots
,\mathfrak{q}_{m}$, such that ${\rm Spec}(A)=\cup
_{i=1}^{m}U(\mathfrak{q}_{i})$. Hence, ${\rm rt}(\mathfrak{m};M)\leq
{\rm max}\{ {\rm rt}(\mathfrak{q}_{i};M)\mid i=1,\ldots ,m\} $ for any
maximal ideal $\mathfrak{m}$ of $A$. \fidemo

\medskip

Using Theorem 2 in \cite{planas2} we deduce the result of Duncan and
O'Carroll in \cite{do}.

\begin{corollari}\label{cucu} {\rm \cite{do}}
Let $A$ be an excellent (or $J-2$) ring and let $N\subseteq M$ be
two finitely generated $A$-modules.
Then there exists an integer $s\geq 1$ such that, for all integers
$n\geq s$ and for all maximal ideals $\mathfrak{m}$ of $A$,
$\mathfrak{m}^{n}M\cap N=\mathfrak{m}^{n-s}(\mathfrak{m}^{s}M\cap N)$.
\end{corollari}

\noindent {\sc Acknowledgement}. This work was partially supported by
the DGES~PB97-0893 grant.

\end{document}